\documentclass[12pt]{amsart}
\usepackage[all]{xy}
\usepackage{amssymb}
\usepackage{amsfonts,amsthm,amsmath,amscd}
\begin{document}{}
\hfuzz=30pt
\newtheorem{theorem}{Theorem}[section]
\newtheorem{lemma}[theorem]{Lemma}
\newtheorem{proposition}[theorem]{Proposition}
\newtheorem{corollary}[theorem]{Corollary}

\theoremstyle{definition}
\newtheorem{definition}[theorem]{Definition}
\newtheorem{example}[theorem]{Example}
\newtheorem{xca}[theorem]{Exercise}

\newtheorem{remark}[theorem]{Remark}

\numberwithin{equation}{section}
\newcommand{\M}{{\text{Mod}}}
\newcommand{\G}{{\textbf{ G}}}
\newcommand{\A}{{\mathcal{A}}}
\newcommand{\B}{{\mathcal{B}}}
\newcommand{\e}{{\epsilon}}
\newcommand{\dd}{{\delta}}
\newcommand{\ta}{{\theta}}
\newcommand{\BB}{{\mathcal{B}_ {\textbf{G}}}}
\newcommand{\SB}{{S \otimes_B S}}
\newcommand{\E}{{\text{End}_{A-\text{cor}}(\Sigma)}}
\newcommand{\Et}{{ \text{End}_{A-\text{cor}}(\Sigma)}}
\newcommand{\Es}{{\text{End}_{B-\text{cor}}(\SB)}}

\title{On a generalization of Grothendieck's theorem}

%Information for first author

\author{Bachuki Mesablishvili}
%Address of record for the research reported here
\address{Razmadze Mathematical Institute, Tbilisi 0193, Republic of Georgia}
\email{bachi@rmi.acnet.ge}
%\thanks{GRDF GEM1-3330-TB-03 grants}

%General info
%\amsclass{18D20}
%\date{ }

%\dedicatory{ }

\begin{abstract} A wide generalization of the classical theorem
of A. Grothendieck asserting that for any faithfully flat
extension of commutative rings, the corresponding relative Picard
group and the Amitsur 1-cohomology  group with values in the
units-functor are isomorphic, is obtained. This implies some known
results that are concerned with extending to non-commutative rings
of Grothendieck's theorem.

\end{abstract}
\maketitle
\section{Introduction}

One of the fundamental results in descent theory is Grothendieck's
theorem (see Corollary 4.6. in \cite{CW}) establishing an
isomorphism between the relative Picard group $\text{Pic}(S/B)$ of
a faithfully flat extension $i: B \to S$ of commutative rings and
the Amitsur 1-cohomology group $H^1(S/B,U)$ of the extension with
values in the units-functor $U$.

Grothendieck's result was generalized in \cite{M} to
non-commutative rings as follows: Let $i: B \to S$ be an extension
of non-commutative rings, let $\text{Inv}_R(S)$ denote the group
of invertible $B$-subbimodules of $S$, and
${\text{Aut}_{A-\text{cor}}(S \!\!\otimes_B \!\!S)}$ the group of
$B$-coring automorphisms of the Sweedler's canonical $B$-coring $S
\!\!\otimes_B \!\!S$. Masouka defined a group homomorphism
$\Gamma:\text{Inv}_R(S) \to {\text{Aut}_{A-\text{cor}}(S
\!\!\otimes_B \!\!S)}$ and showed that if either (a) $S$ is
faithfully flat as a right or left $B$-module, or (b) $B$ is a
direct summand of $S$ as a $B$-bimodule, then $\Gamma$ is an
isomorphism of groups.

This has been further generalized by L. El Kaoutit and J.
G\'{o}mez-Torrecillas \cite{EG}, considering extensions of
non-commutative rings of the form $B \to \text{End}_A(M)$, where
$M$ is a $B$-$A$-bimodule with $M_A$ finitely generated and
projective.

In the present paper, we obtain a more general result that
includes the above results as particular cases.

We refer to \cite{BW} for terminology and general results on
(co)monads, and to \cite{BW1} for a comprehensive introduction to
the theory of corings and comodules.
\section{Preliminaries}
We begin by recalling that a comonad $\G$ on a given category $\B$
is an endofunctor $G : \B \to \B$ equipped with natural
transformations $\e : G \to 1$ and $\dd : G \to G^2$ such that the
diagrams
$$\xymatrix{
G \ar[r]^{\dd} \ar[d]_{\dd} & G^2 \ar[d]^{\dd G} & \text{and}& G &
G^2 \ar[l]_{\e G}
\ar[r]^{G \e} &G\\
G^2 \ar[r]_{G \dd} & G^3 &&& G \ar[u]^{\delta} \ar @{=}[ur]
\ar@{=}[ul] }
$$ are commutative. If $\G =(G, \dd, \e)$ is a comonad on $\B$,
then a $\G$-coalgebra is a pair $(b, \ta _b)$ with $b \in \B$ and
$\ta _b :b \to G(b)$ a morphism in $\B$ for which $\e _b \cdot \ta
_b =1$ and $\dd _b \cdot \ta _b = G(\ta _b ) \cdot \ta _b .$ If
$(b, \ta _b)$ and $(b', \ta _{b'})$ are $\G$-coalgebras, then
their morphism $ f:(b, \ta _b) \to (b', \ta _{b'})$ is a morphism
$f : b \to b'$ of $\B$ for which $\ta _{b'} \cdot f = G(f) \cdot
\ta _b .$

The $\G$-coalgebras and their morphisms form a category $\BB$, the
category of $\G$-coalgebras (or the \emph{Eilenberg-Moore}
category associted to $\G$). There are functors $F_{\G}: \BB \to
\B$ and $U_{\G}: \B \to \BB$, given on objects by $F_{\G}(b, \ta
_b)=b$ and $ U_{\G}(b)= (G(b), \dd _b)$. Moreover,
$F_{\textbf{G}}$ is left adjoint to $U_{\G}$.

Recall also that if $\eta, \epsilon: F \dashv U : \B \to \A$ is an
adjunction (so that $F: \A \to \B$ is a left adjoint of $U : \B
\to \A$ with unit $\eta :1 \to UF$ and counit $\e : FU \to 1$),
then $\G =(G, \e, \dd)$ is a comonad on $\B$, where $G=FU$, $\e :
G=FU \to 1$ and $\dd =F \eta U : G=FU \to FUFU=G^2$, and one has
the comparison functor $K_{\textbf{G}}$ in

$$
\xymatrix@!=2.9pc{
&  \A \ar[dr]^{K_{\textbf{G}}} \ar@<.7ex> [dl]^F & \\
\B \ar@<.7ex> [ur]^U \ar@<-.7ex> [rr]_{U_{\textbf{G}}} && \B
_{\textbf{G}} \ar@<-.7ex> [ll]_{F_{\textbf{G}}}.}
$$
where $K_{\textbf{G}}(a)= (F(a), F(\eta_a ))$ and
$K_{\textbf{G}}(f)=F(f)$. Moreover, $F_{\textbf{G}} \cdot
K_{\textbf{G}} \simeq F$ and  $K_{\textbf{G}} \cdot U \simeq
U_{\textbf{G}}$. One says  that the functor $F$ is
\emph{precomonadic} if $K_{\textbf{G}}$ is full and faithful, and
it is \emph{comonadic} if $K_{\textbf{G}}$ is an equivalence of
categories.

\begin{theorem} \emph{(}Beck, see \cite{BW}\emph{)}  Let $ \eta, \epsilon \!: \!F
\!\dashv \! U \!: \!\B \to \A$ be an adjunction, and let $\G=(FU,
\e, F \eta U )$ be the corresponding comonad on $\B$. Then:

\begin{itemize}

\item[1.] The comparison functor $K_{\bf{G}}: \A \to \B_{\bf {G}}$
has a right adjoint $R_ {\bf {G}}: \B _{\bf {G}} \to \A$ iff for
each $(b, \ta _b ) \in \B_{\bf {G}}$, the pair of morphisms
$(U(\ta _b) , \eta _{U(b)})$ has an equalizer in $\A $ - one then
finds $R_{\bf {G}} (b, \ta _b )$ as the equalizer
\begin{equation}
\xymatrix { R_{\bf {G}} (b, \ta _b )\ar[rr]^{e_{(b, \ta _b )}} &&
U(b) \ar@{->}@<0.5ex>[rr]^{U(\ta _b)} \ar@ {->}@<-0.5ex>
[rr]_{\eta _{U(b)}} && UFU(b).}
\end{equation}

\item[2.] Assuming the existence of $R_ {\bf {G}}$, $K_ {\bf {G}}$
is an equivalence of categories (in other words, $F$ is comonadic)
iff the functor $F$ is conservative(=isomorphism-reflecting) and
preserves (or equivalently, preserves and reflects) the equalizer
$(2.1)$ for each $(b, \ta _b ) \in \B_{\bf {G}}$.
\end{itemize}
\end{theorem}

\bigskip

Let $i : B \to S$ be an arbitrary extension of (non-commutative)
rings, $\A$ be the category ${_B \M}$ of left $B$-modules, $\B$ be
the category ${_S \M}$ of left $S$-modules, $$F_S =S \otimes _B :
{_B \M} \to {_S \M}$$ and $$U_S  : {_S \M} \to {_B \M}$$ be the
restriction-of-scalars functor. It is well known that $F_S$ is
left adjoint to $U_S$ and that the unit $\eta$ of this adjunction
is given by $$\eta_X : X \to S \otimes _B X, \, \eta_X (x)=1
\otimes_B x .$$

It is also well known that the Eilenberg-Moore category ${({_S
\M})\!\!_{\G}}$ of $\G$-coalgebras,$\G$ being the comonad on ${_S
\M}$ associated to the adjunction $F_S \dashv U_S$, is equivalent
to the category ${^{\SB}({_S \M})}$ of left comodules over the
Sweedler canonical $B$-coring $\SB$ corresponding  to the ring
extension $i$, by an equivalence which identifies the comparison
functor $K_{\G} : {_B \M} \to {{({_S \M})_{\G}}}$ with the functor
$$K_S : {_B \M} \longrightarrow {^{\SB}({_S \M})}, \,\, K_S (X)=(S \otimes_B X,
\ta _{S \otimes_B X}),$$ where $\ta _{S \otimes_B X}= S \otimes_B
\eta_X$ for all $X \in {_B \M}$. (Note that a left $\SB$-comodule
is a pair $(Y, \ta_Y)$ with $Y \in {_S \M}$ and $\ta_Y : Y \to S
\otimes_B Y$ a left $A$-module morphism for which the diagrams
$$
\xymatrix { Y \ar@{=}[rd] \ar[r]^-{\ta_Y} & S \otimes_B Y
\ar[d]^{\alpha_Y} &\text{and}& Y \ar[d]_{\ta_Y}\ar[r]^-{\ta_Y} & S
\otimes_B Y \ar[d]^{S
\otimes_B \eta_Y}\\
&Y\,\, && S \otimes_B Y \ar[r]_-{S \otimes_B \ta_Y} &S \otimes_B S
\otimes_B Y,}
$$ where $\alpha_Y$ denotes the left $S$-module structure on $Y$, are
commutative.) So, to say that the functor $F_S =S \otimes_B -$ is
comonadic is to say that the functor $K_S$ is an equivalence of
categories. Applying Beck's theorem and using that ${_B \M}$ has
all equalizers, we get:

\begin{theorem} The functor $F_S =S \otimes_B - : {_B \emph{\M}}
\to {_S \emph{\M}}$ is comonadic if and only if
\begin{itemize}

\item[(i)]the functor $F_S$ is conservative, or equivalently, the
ring extension $i: B\to S$ is a \emph{pure} morphism of right
$B$-modules;

\item[(ii)] for any $(Y, \ta_Y) \in {^{\SB}({_S \emph{\M}})}$,
$F_S$ preserves the equalizer
\begin{equation}
\xymatrix { R_S(Y, \ta_Y) \ar[rr]^-{e_{(Y, \ta_Y)}} && Y
\ar@{->}@<0.5ex>[r]^-{\eta_Y} \ar@ {->}@<-0.5ex> [r]_-{\ta_Y} &
S\otimes_B Y,}\end{equation} where $R_S :{^{\SB}({_S \emph{\M}})}
\to {_B \emph{\M}}$ is the right adjoint of the comparison functor
$K_S : {_B \emph{\M}} \to {^{\SB}({_S \emph{\M}})}$.
\end{itemize}
\end{theorem}

\bigskip
\bigskip
Let $A$ be a ring and $\Sigma$ be an $\A$-coring. Let us write
$\text{End}_{A-\text{cor}}(\Sigma)$ (resp.
$\text{Aut}_{A-\text{cor}}(\Sigma)$) for the monoid (resp. group)
of $A$-coring endomorphisms (resp. automorphisms) of $\Sigma$.
Recall that any $g \in \E$ induces functors: $${_g(-)} :
{^{\Sigma}({_A \M}}) \to {^{\Sigma}({_A \M}}),$$ defined by
${_g(Y, \ta_Y)} =(Y, (g \otimes_A 1)\circ \ta_Y)$, and $$(-)_g :
\M ^{\Sigma} \to \M ^{\Sigma}$$ defined by $(Y',
\ta_{Y'})_g=(Y',(1 \otimes_A g) \circ \ta_{Y'} )$.

\bigskip
\bigskip

It is easy to see that the left $S$-module $S$ is a left
$(\SB)$-comodule with left coaction $${_S \ta} : S \to \SB, \,\,\,
s \to s \otimes_B 1,$$ and that ${_g(S, {_S \ta})=(S, g\circ {_S
\ta}})$. Symmetrically, the right $S$-module $S$ is a right
$\SB$-comodule with the right action $$\ta_S : S \to S\otimes_B S,
\,\, s \to 1 \otimes_B s,$$ and that $(S, \ta_S)_g= (S, g \circ
\ta_S)$.

\bigskip
\bigskip
For a given injective homomorphism $i: B \to S$ of rings, let

\begin{itemize}

\item $I_B (S)$ denote the monoid of all $B$-subbimodules of $S$,
the multiplication being given by $$IJ=\{\sum_{k\in K} i_k \cdot
j_k, \,\,\, I,J \in I_B (S),\,\, i_k \in I, \,\, j_k \in J\, , \,
\text{and}\, K \, \text{is a finite set}\};$$

\item $I^l _B (S)$ (resp. $I^r _B (S)$) denote the submonoid of
$I_B (S)$ consisting of those $I \in I_B (S)$ for which the map
$$\textbf{m}^l_I : S \otimes_B I \to S, \,\, m \otimes _B i \to mi,
$$$$\,\,\,(\text{resp.} \,\,\textbf{m}^r_I: I \otimes_B S \to S,\,\, i
\otimes_B m \to im)$$ is an isomorphism;

\item $J(g)=\{s \in S \,\,|\,\, g(s\otimes_B 1)=1 \otimes_B s\}
\,\,\, \text{for}\,\, g \in \Es$ and let ${i_g} : J(g) \to S$ be
the canonical embedding;

\item $J'(g)=\{s \in S \,\,|\,\, s\otimes_B 1=g(1 \otimes_B s)\}
\,\,\, \text{for}\,\, g \in \Es$ and let ${i'_g} : J'(g) \to S$ be
the canonical embedding.
\end{itemize}

It is clear that $J(g), \, J'(g) \in I_B(S)$ for all $ g \in \Es$.

\bigskip
\bigskip
The following result is verified directly:
\begin{proposition} For any $ g \in \Es$, $R_S ({_g(S, {_S \ta})})
\simeq J(g)$.
\end{proposition}

\section{Main Results}

In this section we present our main results.
\bigskip

We begin with
\begin{proposition}For any $g \in \Es$, the following conditions
are equivalent:
\begin{itemize}
\item [(i)] $ J(g) \in I^l_B (S);$

\item [(ii)] the ${_g(S, {_S \ta})}$-component of the counit
$\varepsilon : K_S R_S \to 1$ of the adjunction $K_S \dashv R_S$
is an isomorphism;

\item [(iii)]the functor $S \otimes_B - : {_B \emph{\M}} \to {_S
\emph{\M}}$ preserves the equalizer \begin{equation}\xymatrix {
J(g)  \ar[r]^-{{i_g}} & S \ar@{->}@<0.5ex>[r]^-{\eta_S} \ar@
{->}@<-0.5ex> [r]_-{g \circ {_S\ta} } & S\otimes_B
S;}\end{equation}

\item [(iv)]the morphism $S \otimes_B \!i_g : S \otimes_B J(g) \to
S \otimes_B S$ is a monomorphism.
\end{itemize}
\end{proposition}

\begin{proof}It is well known (see, for example, \cite{BW} ) that, for
any $(Y, \ta_Y) \in {^{\SB}({_S \M})}$, the diagram $$\xymatrix {
Y \ar[r]^-{\ta_Y} & S \otimes_B Y \ar@{->}@<0.5ex>[r]^-{S
\otimes_B \ta_Y} \ar@ {->}@<-0.5ex> [r]_-{S \otimes_B \eta_Y} &
S\otimes_B S \otimes_B Y}$$ is an equalizer and that the $(Y,
\ta_Y)$-component $\varepsilon_{(Y, \ta_Y)}$ of $\varepsilon$
appears as the unique factorization of the morphism $S \otimes_B
e_{(Y, \ta_Y)}$ through the morphism $\ta_Y$:
\begin{equation}\xymatrix {
S \otimes_B \!R_S(Y, \ta_Y) \ar[rd]^-{S \otimes_B e_{(Y, \ta_Y)}}
\ar@{.>}[d]_{\varepsilon_{(Y, \ta_Y)}} &&\\
Y \ar[r]_{\ta_Y}& S \otimes_B Y \ar@{->}@<0.5ex>[rr]^-{S \otimes_B
\ta_Y}\ar@ {->}@<-0.5ex> [rr]_-{S \otimes_B \eta_Y}&&S\otimes_B S
\otimes_B Y.}\end{equation} Since $\alpha_Y \cdot \ta_Y=1$,
$\varepsilon_{(Y, \ta_Y)}=\alpha_Y \cdot (S \otimes_B e_{(Y,
\ta_Y)}).$ In particular, when  $(Y, \ta_Y)={_g(S, {_S \ta})}$ we
get that $\varepsilon_{{_g(S, {_S \ta})}}=\textbf{m}^l_{J(g)}$. So
$\text{(i)}$ and $\text{(ii)}$ are equivalent.

Since the row of the diagram (3.2) is an equalizer, it follows
that the morphism $S \!\otimes_B \!e_{(Y, \ta_Y)}$ is an equalizer
of the pair of morphisms $(S \otimes_B \ta_Y, S \otimes_B \eta_Y)$
iff $\varepsilon_{(Y, \ta_Y)}$ is an isomorphism. In other words,
the functor $S \otimes_B- $ preserves the equalizer (2.2) iff
$\varepsilon_{(Y, \ta_Y)}$ is an isomorphism. As a special case we
then have that $\text{(ii)}$ is equivalent to $\text{(iii)}$.

Finally, since the category $_B \M$ is abelian (and hence coexact
in the sense of Barr \cite{BW}), and since $i_g$ is the equalizer
of the $(S \otimes_B-)$-split pair of morphisms $(_S \ta,
\eta_S)$, it follows from the proof of Duskin's theorem (see, for
example, \cite{BW}) that the functor $S \otimes_B-$ preserves the
equalizer (3.1) iff the morphism $S \otimes_B i_g$ is a
monomorphism. So $\text{(iii)}$ and $\text{(iv)}$ are also
equivalent. This completes the proof.
\end{proof}

It is shown in \cite{EG} that assigning to each $I \in I^l_B (S)$
(resp. $I \in I^r_B (S)$) the composite $\Gamma=(1 \otimes_B
\textbf{m}^r_I) \circ ((\textbf{m}^l_I)^{-1} \otimes_B 1)$ (resp.
$\Gamma'=(\textbf{m}^l_I \otimes_B 1) \circ (1 \otimes_B
(\textbf{m}^r_I)^{-1})$) yields an (anti-)homomorphism of monoids
$\Gamma: I^l_B (S) \to \Es$ (resp. $\Gamma': I^l_B (S) \to \Es$).

We shall need the following easy consequence of Lemma 2.7 of
\cite{M}:
\begin{proposition}Assume that $i: B \to S$ is such
that any embedding $I \hookrightarrow J$ of $B$-subbimodules of
$S$ is an isomorphism whenever its image under the functor $S
\otimes_B -$ is such. Then $\Gamma: I^l_B (S) \to \Es$ is an
isomorphism of monoids whose inverse is the map $ g \to J(g)$,
provided that $J(g) \in I^l_B (S)$ for all $g \in \Es$.
\end{proposition}

Putting Propositions 3.1 and 3.2 together, we get:

\begin{theorem}Let $i: B \to S$ be as in Proposition 3.2. Then
$\Gamma: I^l_B (S) \to \Es$ is an isomorphism of monoids if and
only if, for any $g \in \Es$, the equivalent conditions of
Proposition 3.1 hold.
\end{theorem}

\bigskip

\begin{proposition}If the functor $S \otimes_B - : {_B \emph{\M}} \to
{_S \emph{\M}}$ is comonadic, then $J(g) \in \emph{I}^l _B (S)$
for all $g \in \Es$.
\end{proposition}

\begin{proof} Consider the left $(\SB)$-comodule $(S, {_S\ta})$.
According to Proposition 2.3 and Theorem 2.1, the pair $(J(g),
{i_g} : J(g) \to S)$ appears as the equalizer
$$\xymatrix { J(g)  \ar[r]^-{{i_g}} & S
\ar@{->}@<0.5ex>[r]^-{\eta_S} \ar@ {->}@<-0.5ex> [r]_-{g \circ
{_S\ta} } & S\otimes_B S,}
$$ and since the functor $S\otimes_B -$ is assumed to be comonadic,
it preserves the equalizer $(2.2)$ for all $(Y, \ta_Y) \in
{^{\SB}({_S \M})}$ and in particular considering $(S, {_S\ta}) \in
{^{\SB}({_S \M})}$, we see that

$$\xymatrix { S \otimes_B J(g) \ar[rr]^-{S \otimes_B {i_g}} && S
\otimes_B S \ar@{->}@<0.5ex>[rr]^-{S \otimes_B \eta_S} \ar@
{->}@<-0.5ex> [rr]_-{S \otimes_B (g \circ {_S\ta})} && S\otimes_B
S \otimes_B S}
$$ is an equalizer diagram. It now follows from Proposition 3.1
that $J(g) \in I^l _B (S)$.
\end{proof}

Recalling that any comonadic functor is conservative, and putting
Theorem 3.3 and  Proposition 3.4 together, we obtain:

\begin{theorem}If the functor $S \otimes_B - : {_B \emph{\M}} \to {_S \emph{\M}}$
is comonadic, then $\Gamma : I^l _B (S) \to \Es$ is an isomorphism
of monoids.
\end{theorem}

There is of course a dual result.

\begin{theorem}If the functor $- \otimes_B S : { \emph{\M}_B} \to
{\emph{\M}_S}$ is comonadic, then $\Gamma' : I^r _B (S) \to \Es$
is an anti-isomorphism of monoids.
\end{theorem}

It is known (see \cite{M}) that the monoid morphism $$\Gamma :
\text{I}^l _B (S) \to \Es$$ restricts to a group morphism $$
\text{Inv}_B (S) \to \text{Aut}_{B-{cor}}(\SB),$$ which is still
denoted by $\Gamma $.

\begin{theorem} If either

\begin{itemize}

\item[(i)] the functor $S \otimes_B - : {_B \emph{\M}} \to {_S
\emph{\M}}$, or

\item[(ii)] the functor $- \otimes_B S : { \emph{\M}_B} \to
{\emph{\M}_S}$
\end{itemize}
is comonadic, then $\Gamma : \emph{Inv} _B (S) \to
\text{Aut}_{B-\text{cor}}(\SB)$ is an isomorphism of groups.
\end{theorem}

\begin{proof}The same argument as in \cite{EG} shows that if either
$\Gamma : I^l _B (S) \to \Es$ or $\Gamma' : I^r _B (S) \to \Es$ is
an isomorphism, then the group homomorphism $\Gamma$ is an
isomorphism. Theorems 3.5 and 3.6 now complete the proof.
\end{proof}

As a special case of this theorem, we obtain the following result
of Masuoka (see \cite{M}):

\begin{theorem}If either
\begin{itemize}
\item[(i)] ${_B S}$ is faithfully flat, or

\item[(ii)] $B$ is a direct summand of $S$ as a $B$-bimodule,

\end{itemize}
then $\Gamma : I^l _B (S) \to \Es$ is an isomorphism of monoids.
\end{theorem}

\begin{proof} In both cases, the functor $S \otimes_B - :
{_B \M} \to {_S \M}$ is comonadic. Indeed, to say that ${_B S}$ is
faithfully flat is to say that the functor $S \otimes_B -$ is
conservative and it preserves all equalizers. Thus, according to
Beck's theorem, this functor is comonadic.

Now, if $B$ is a direct summand of $S$ as a $B$-bimodule, it is
not hard to see that the unit of the adjunction $F_S=S \otimes_B -
\dashv U_S $ is a split monomorphism and it follows from Theorem
2.2 of \cite{JT}  that the functor $F_S$ is comonadic. Theorem 3.7
now completes the proof.
\end{proof}

Dually we have:
\begin{theorem}If either
\begin{itemize}
\item[(i)] ${S_B}$ is faithfully flat, or

\item[(ii)] $B$ is a direct summand of $S$ as a $B$-bimodule,

\end{itemize}
then $\Gamma ' : I^r _B (S) \to \Es$ is an anti-isomorphism of
monoids.
\end{theorem}

\bigskip

\begin{theorem} If either
\begin{itemize}
\item[(i)] ${_B S}$ or $S_B$ is faithfully flat, or

\item[(ii)] $B$ is a direct summand of $S$ as a $B$-bimodule,

\end{itemize}
then $\Gamma : \emph{Inv} _B (S) \to
\text{Aut}_{B-\text{cor}}(\SB)$ is an isomorphism of groups.
\end{theorem}

\begin{proof} The argument here is the same as in the proof of
Theorem 3.7.
\end{proof}

We now consider the following situation: Let $A$ and $B$ be rings,
$M$ a $(B,A)$-bimodule with $M_A$ finitely generated and
projective, $S=\text{End}_A (M)$ the ring of right
$A$-endomorphisms of $M_A$, and $\Sigma =M^* \otimes_B M$ the
comatrix $A$-coring corresponding to ${_B M_A}$ (for the notion of
comatrix coring see \cite{EG2}). When ${_B M_A}$ is faithful, in
the sense that the canonical morphism
$$i : B \to S, \,\, s \longrightarrow [m \to sm]$$ is injective,
one has a map
$$\Gamma_0  : I^l _B (S) \to \Et$$ of sets defining
$\Gamma_0 ^l(I), \,\, I \in I^l _B (S),$ to be the endomorphism
$$m^* \otimes_B m \longrightarrow \sum_i m^* x_i \otimes_B y_i
m,$$ where $(\textbf{m}^l_I)^{-1}(1)=\sum_i x_i \otimes_B y_i \in
I^l _B (S).$

\begin{theorem}Suppose that ${_B M_A}$ is such that the functor
$$S \otimes_B - :{_B \emph{\M}} \to {_S \emph{\M}}$$ is comonadic.
Then the map $$\Gamma_0  : I^l _B (S) \to \Et$$ is in fact an
isomorphism of \emph{monoids}.
\end{theorem}

\begin{proof} First of all, the morphism $i: B \to S$ is
injective (or equivalently, the bimodule ${_B M_A}$ is faithful),
since the functor $S \otimes_B -$ is assumed to be comonadic.
Next, it is proved in \cite{EG} that the assignment
$$g \longrightarrow \widehat{g}=(\xi \otimes_B \xi) \circ (M
\otimes_A g \otimes_A M^*) \circ (\xi^{-1} \otimes_B \xi^{-1}),$$
where $\xi: M\otimes_A M^* \to S=\text{End}_A (M)$ is the
canonical isomorphism, yields an injective morphism of monoids
$$\widehat{(-)}: \Et \to \Es.$$ And the same argument as in the
proof of Proposition 2.6 of \cite{EG} shows that the following
diagram of \emph{sets}

$$\xymatrix {
I^l _B (S) \ar[d]^{\Gamma } \ar[r]^{\Gamma_0 } &
\Et \ar[dl]^{\hat{(-)}}\\
\Es }
$$ is commutative. Now, since the functor $S \otimes_B -$ is
assumed to be comonadic, it follows from Theorem 3.5 that $\Gamma
$ is an isomorphism of monoids and hence the monoid morphism
$\widehat{(-)}$, being injective, is also an isomorphism.
Commutativity of the diagram then gives that $\Gamma _0$ is an
isomorphism of monoids.
\end{proof}

Dually, one can define a map  $$\Gamma'_0  : I^r _B (S) \to \Et$$
that sends $I \in I^r _B (S)$ to the endomorphism
$$ m^* \otimes_B m \longrightarrow \sum_i m^*y_i \otimes_B x_i
m$$ of the $A$-coring $\Et$, where
$(\textbf{m}^r_I)^{-1}(1)=\sum_i y_i \otimes_b x_i \in I \otimes_B
S.$

\begin{theorem}Suppose that ${_B M_A}$ is such that the functor
$$- \otimes_B S :{ \emph{\M}_B} \to  \emph{\M}_S$$ is comonadic. Then
$$\Gamma'_0  : I^r _B (S) \to \Et$$ is an anti-isomorphism of
\emph{monoids}.
\end{theorem}

\bigskip

It is not hard to check that the map $$\Gamma_0  : I^l _B (S) \to
\Et$$ of sets restricts to a map $$\text{Inv}_B (S) \to
\text{Aut}_{A-\text{cor}}(\Sigma)$$ which we still call
$\Gamma_0.$ As in \cite{EG} , it follows from Theorems 3.11 and
3.12 that

\begin{theorem} If either
\begin{itemize}

\item [(i)] the functor $S \otimes_B -$, or

\item [(ii)] the functor $- \otimes_B S$
\end{itemize}
is comonadic, then the map $$\Gamma_0 : \emph{Inv}_B (S) \to
\text{Aut}_{A-\text{cor}}(\Sigma)$$ is actually an isomorphism of
groups.
\end{theorem}

It is shown in \cite{Me} that the functor $S \!\otimes_B -:{_B
\M}\to {_S \M}$ (resp. $- \otimes_B \!S: \M_B \to \M_S$) is
comonadic iff the functor $M \!\otimes_B -: {_B \M}\to {_A \M}$
(resp. $- \otimes_B \!M^*: \M_B \to \M_A$) is. So we have:

\begin{theorem} If either
\begin{itemize}

\item [(i)] the functor $M \!\otimes_B -$, or

\item [(ii)] the functor $- \otimes_B \!M^*$
\end{itemize}
is comonadic, then the map $$\Gamma_0 : \emph{Inv}_B (S) \to
\text{Aut}_{A-\text{cor}}(\Sigma)$$ is an isomorphism of groups.
\end{theorem}

From the last theorem one obtains the following result of L. El
Kaoutit and J. G\'{o}mez-Torrecillas (see Theorem 2.5 in
\cite{EG}):

\begin{theorem}If

\begin{itemize}

\item[(i)] ${_B M}$ is faithfully flat, or

\item[(ii)] ${M^*_B}$ is faithfully flat, or

\item[(i)] ${_B M_A}$ is a separable bimodule,
\end{itemize}
then $$\Gamma_0 : \text{Inv}_B (S) \to
\text{Aut}_{A-\text{cor}}(\Sigma)$$ is an isomorphism of groups.
\end{theorem}

\begin{proof} $(\text{i})$ and $(\text{ii}).$ To say that ${_B M}$ (resp.
$M^*_B$) is faithfully flat is to say that the functor $M
\!\otimes_B -: {_B \M}\to {_A \M}$ (resp. $- \otimes_B \!M^*: \M_B
\to \M_A$) is conservative and preserves all equalizers. Then the
functor $S\otimes_B - : {_B \M} \to {_B \M}$ (resp. $S- \otimes_B
S : \M_B \to \M_S$) is comonadic by a simple application of the
Beck theorem. Applying the previous theorem, we see that $\Gamma$
is an isomorphism of groups.

$(\text{iii}).$ If ${_B M}_A$ is a separable bimodule, then the
ring extension $i: B \to S$ splits (see, for example, \cite{S}),
i.e. $B$ is a direct summand of $S$ as a $B$-bimodule. But we have
already seen (see the proof of Theorem 3.8) that in this case, the
functor $S\otimes_B - : \M_B \to \M_S$ is comonadic, and Theorem
3.13 shows that $\Gamma$ is an isomorphism of groups.
\end{proof}

\bibliographystyle{amsplain}

\end{document}